\documentclass{article}
\usepackage{amssymb}
\usepackage{amsmath}

\input{tcilatex}
\begin{document}

\bigskip

\begin{center}
\noindent {\LARGE Control of a finite dam when the input process is either
spectrally positive L\'{e}vy or spectrally positive L\'{e}vy reflected at
its infimum}

{\LARGE \bigskip }

Mohamed Abdel-Hameed

Department of Statistics

College of Business and Economics

UAE University
\end{center}

\bigskip

\bigskip {\LARGE Abstract}

\smallskip\ Bae \textit{et al}. [6] consider the problem of optimal control
of a finite dam using $P_{\lambda ,\tau }^{M}$ policies, assuming that the
input process is a compound Poisson process with a negative drift. Lam and
Lou [8] treat the case where the input is a Wiener process with a reflecting
boundary at its infimum, with drift term $\mu \geq 0$, using the long-run
average and total discounted cost criteria. Attia [4] obtains results
similar to those of Lam and Lou, through simpler and more direct methods.
Zuckermann [12] considers $P_{\lambda ,0}^{M}$ policies when the input
process is a is Wiener process with drift term \ $\mu \geq 0$. The
techniques used by the above mentioned authors\ involve solving systems of
differential or integral equations. In this paper we use the theory and
methods of scale functions of L\'{e}vy processes to unify and extend the
results of these authors.



\bigskip

\medskip \noindent {\LARGE 1}. {\LARGE Introduction and summary} \bigskip 
\newline
\indent\ Suppose that a dam has capacity $V$. Its water input $%
I=(I_{t},t\geq 0)$, is assumed to be a L\'{e}vy\ process with drift $\mu $,
variance $\sigma ^{2}$, and the water is released at one of two rates $0$ or 
$M$ \ units per unit of time. We consider $P_{\lambda ,\tau }^{M}\ $policies
in which the water release rate is assumed to be zero until the water
crosses level $\lambda $, $(0<\lambda <V)$, when the water is released at
rate $M\ \ $until it reaches level $\tau $, $(0\leq \tau <\lambda )$. Once
level $\tau \ $is reached, the release rate remains zero until level $%
\lambda \ $is reached again, and the cycle is repeated. We deal with the
cases where the input process is spectrally positive L\'{e}vy, and
spectrally positive L\'{e}vy reflected at its infimum. In both cases the
content process is a delayed regenerative process with regeneration points
being the times of successive visits to state $\tau $. \ During a given
cycle, the dam's water content is a L\'{e}vy process with coefficients $\mu $
and $\sigma ^{2}$, and it remains so until it crosses level $\lambda $; from
then until it drops to level $\tau \ $again the content level behaves like a
L\'{e}vy process\ reflected at V with coefficients $\mu ^{\ast }=\mu -M$,$\
\sigma ^{2}$, denoted by $I^{\ast }=(I_{t}^{\ast },t\geq 0)$. At any time,
the release rate can be increased from $0$ to $M$ with a starting cost$%
\;K_{1}M$, or decreased from$\;M$\ to zero with a closing cost $K_{2}M$.\
Moreover, for each unit of output, a reward $R\;$is received.\ Furthermore,
there is a penalty cost which accrues at a rate $f$ , where $f$ \ a bounded
measurable function. \ For the first case we extend the results of Zuckerman
[12] who assumed that $\tau =0$ and $f=0$. Our results in the second case
extend the results of Lam and Lou [8] and Attia [4], where they assumed that
the input process is a Wiener process reflected at its infimum. They also
extend those of Bae \textit{et al. }[6], who consider the case where the
water input is a compound Poisson process with negative drift. Lee and Ahn
[9] consider the long-run average cost case, for the $P_{\lambda ,0}^{M}\ $%
policy, when the water input is a compound Poisson process. Abdel-Hameed [1]
treats the case where the water input is a compound Poisson process with a
positive drift. He obtains the total discounted as well as the long-run
average costs. Bae \textit{et al. }[5\ ] consider the $P_{\lambda ,0}^{M}\ $%
policy in assessing the workload of an M/G/1 queuing system. The techniques
used by in [12], [8], [4], and [6] \ involve solving systems of differential
or integral equations. In this paper we use the theory and methods of scale
functions of L\'{e}vy processes, an approach not uses by researchers in this
area before.

In Section 2 we define the input processes and discuss their properties. In
Section 3 we obtain formulas needed for computing the cost functionals. In
Section\ 4,\ we discuss the cost functionals using the total discounted as
well as the long-run average cost cases. In section 5 we discuss the special
cases where the input process is a Gaussian process, a Gaussian process
reflected at its infimum and a spectrally positive L\'{e}vy process of
bounded variation.

\bigskip \smallskip

\bigskip \bigskip \bigskip \noindent {\LARGE 2.~Spectrally positive L\'{e}vy
processes and scale functions}

\bigskip In this section we give some basic definitions; describe spectrally
positive L\'{e}vy processes and discuss some of their characteristics.\ The
reader is referred to [7] for a more detailed discussion of the definitions
and results mentioned in this section.

\medskip

\smallskip For any process $Y=\{Y_{t},t\geq 0\}\ $with state space $E$,$\;$%
any Borel set \ $A$ $\subset E\ \ $and any functional $f$, $E_{y}(f)\;$%
denotes the expectation of\ $f$\ \ conditional on $Y_{0}=y$,$\;P_{y}(A)\;$%
denotes the corresponding probability measure and $\mathbf{I}_{A}(\;)\;$is
the indicator function of the set $A$. In the sequel we will write
indifferently $P_{0}\ $or $P\ \ $and $E\ _{0}$\ or $E$. Throughout, we let $%
R=(-\infty ,\infty )$, $R_{+}=[0,\infty ),$\ $N=\{1,2,...\}$ and $%
N_{+}=\{0,1,...\}$. For $x,y\in R$, we define $x\vee y=x\max y$ $\ $and $%
x\wedge y=x\min y$. For every $t\geq 0$, we define $\underline{Y_{t}}=%
\underset{0\leq s\leq t}{\inf }(Y_{s}$,$\wedge 0),\ \overset{-}{Y}_{t}=%
\underset{0\leq s\leq t}{\sup }(Y_{s},\vee 0)$.

We will use the term "increasing" to mean "non-decreasing" throughout this
paper.

\medskip

\textbf{Definition 1.} A L\'{e}vy process \ $L=$ $\{L_{t},t\geq 0\}$ with$\ $
state space $R\ $is said to be spectrally positive L\'{e}vy process, if it
has no negative jumps.\newline

\smallskip

It follows that, for each $\theta \in R_{+},x\in R$,

\begin{equation*}
E[e^{-\theta \ L\ _{t}}]=e^{t\phi (\theta )},
\end{equation*}%
where

\begin{equation*}
\phi (\theta )=-a\theta +\frac{\theta ^{2}\sigma ^{2}}{2}-\int_{0}^{\infty
}(1-e^{-\theta x}-\theta x\mathbf{I}_{\{x<1\}})\upsilon (dx)\text{. \ \ }%
(2.1)
\end{equation*}%
The terms $a\in R$, $\sigma ^{2}\in R_{+}\ $are\ the drift and variance of
the spectrally positive L\'{e}vy process, respectively. The L\'{e}vy measure 
$\upsilon $ is a positive measure on $(0,\infty )$ satisfying $%
\int_{0}^{\infty }(x^{2}\wedge 1)\upsilon (dx)<\infty $.

\smallskip

The function $\phi \ $is known as the L\'{e}vy exponent, and it is strictly
convex and tends to infinity as $\theta \ $tends to infinity. For $\alpha
\in R_{+}$, we define

\begin{equation*}
\eta (\alpha )=\sup \{\theta :\phi (\theta )=\alpha \}\ \ (2.2)\text{,}
\end{equation*}%
the largest root of the equation $\phi (\theta )=\alpha $. It is seen that
this equation has at most two roots, one of which is the zero root. Note
that,$\ E(L_{1\ })=$ $\int_{1}^{\infty }x\upsilon (dx)+\mu $.$\ $%
Furthermore, $\underset{t\rightarrow \infty }{\lim }L_{t\ }\ =\infty \ $if
and only if $E(L_{1\ })>0,\ $and $\underset{t\rightarrow \infty }{\lim }%
L_{t\ }\ =-\infty \ $if and only if $E(L_{1\ })<0$. Also, if$\ E(L_{1\
})=0,\ $then $\underset{t\rightarrow \infty }{\lim }L_{t\ }$\ does not
exist. Furthermore, $\eta (0)>0,$ if and only if $E(L_{1\ })>0.$

\smallskip

\smallskip An important case is when the process $L\ $is of bounded
variations, i.e., $\sigma ^{2}=0$ and $\int_{0}^{\infty }(x\wedge 1)\upsilon
(dx)<\infty $. Let

\begin{equation*}
\zeta =-a+\int_{0}^{1}x\upsilon (dx)\text{.}
\end{equation*}%
In this case we can write

\begin{equation*}
\phi (\theta )=\zeta \theta -\int_{0}^{\infty }(1-e^{-\theta x})\upsilon (dx)%
\text{, }(2.3)\text{ }
\end{equation*}%
where necessarily $\zeta $ is strictly positive. \newline

\smallskip

\textbf{Definition 2.} A L\'{e}vy process \ is said to be spectrally
negative if it has no positive jumps. \ \ 

\ \ \ \ \ \ \ \ \ \ \ \ \ \ \ \ \ \ \ \ \ \ \ \ \ \ \ 

\textbf{Definition 3. }For any spectrally positive\ \ L\'{e}vy input
process\ $L$, we let $\overset{\wedge }{L}$\ $\ =-\ L$ .

It is clear that $L$ is \ spectrally positive if and only if the process $%
\overset{\wedge }{L}\ \ $is spectrally negative.

\medskip

We now introduce tools, which will be central in the rest of this paper.

\smallskip

\textbf{Definition 4.} For any spectrally positive L\'{e}vy\ process with L%
\'{e}vy exponent $\phi \ $and for $\alpha \geq 0$, the $\alpha -$\textit{\
scale function} $W^{\alpha }:R\twoheadrightarrow R_{+}$, $W^{\alpha }(x)=0\ $%
for every $x\,<\,0,$ and on $[0,\infty )$ it is defined as the unique right
continuous increasing function such that

\begin{equation*}
\int_{0}^{\infty }e^{-\beta x}W^{(\alpha )}(x)dx=\frac{1}{\phi (\beta
)-\alpha }\text{, }\beta >\eta (\alpha )\ \ (2.4)
\end{equation*}%
\smallskip We will denote $W^{0\ \text{ }}$by $W$ throughout. For $\alpha
\geq 0$, we have (see (8.24) of [ 7])

\begin{equation*}
W^{(\alpha )}(x)=\sum_{k=0}^{\infty }\alpha ^{k}W^{\ast (k+1)}(x)\text{,}\ \
\ \ (2.5)
\end{equation*}%
where $W^{\ast (k)}\ $is the kth convolution of $W\ \ $with itself. \ \ \ \ 

\ \ \ \ \ \ \ \ \ \ \ \ \ \ \ \ \ \ \ \ \ \ \ \ \ 

It follows that $W^{(\alpha )}(0+)=0\ $if and only if the process $L$\ is of
unbounded variation. Furthermore, $W^{(\alpha )}$ is right and left
differentiable on $(0,\infty ).\ $By $W_{+}^{(\alpha )^{\prime }}(x)$, we
will denote the right derivative of $W^{(\alpha )}\ $in $x$. \ \ 

\ \ \ \ \ \ \ \ \ \ \ \ \ \ \ \ \ \ \ \ \ \ \ \ \ \ \ 

The adjoint $\alpha -$\ scale function associated\textit{\ with }$W^{(\alpha
)}\ ($denoted$\ $by $Z^{(\alpha )}$) is defined as follows:

\ \ \ \ \ \ \ \ \ \ \ \ \ \ \ \ \ \ \ \ \ \ \ \ \ \ \ \ \ 

\textbf{Definition 5.} For $\alpha \geq 0$, the \textit{adjoint }$\alpha -$%
\textit{\ scale }$Z^{(\alpha )}:R_{{}}\rightarrow \lbrack 1,\infty )\ $is
defined as

\begin{equation*}
Z^{(\alpha )}(x)=1+\alpha \int_{0}^{x}W^{\alpha }(y)dy.\ \ \ (2.6)
\end{equation*}

\smallskip

\smallskip It follows that as $x\rightarrow \infty $,\ for $\alpha >0$, $%
W^{(\alpha )}(x)\sim \frac{e^{\eta (\alpha )x}}{\phi ^{^{\prime }}(\eta
(\alpha )\ )}$ and $\frac{Z^{\alpha }(x)}{W^{(\alpha )}(x)}\sim \frac{\alpha 
}{\eta (\alpha )}$.

\medskip

\bigskip {\LARGE 3.} {\LARGE Basic results}

\bigskip \noindent For each $t\in R_{+}$, let $Z_{t}\ $ be the dam content
at time $t\ $,$\ Z=\{Z_{t},t\in R_{+}\}$.\ We define the following sequence
of stopping times : 
\begin{eqnarray*}
\overset{\wedge }{T}_{0} &=&\inf \{t\geq 0:Z_{t}\geq \lambda \},\;\;\;\;\;%
\overset{\ast }{T}_{0}=\inf \{t\geq \overset{\wedge }{T}_{0}:Z_{t}=\tau \},
\\[0.01in]
\overset{\wedge }{T}_{n} &=&\inf \{t\geq \overset{\ast }{T}_{n-1}:Z_{t}\geq
\lambda \},\;\overset{\ast }{T}_{n}=\inf \{t\geq \overset{\wedge }{T}%
_{n}:Z_{t}=\tau \}\text{, }n=1,2,...\;(3.1)
\end{eqnarray*}%
It follows that the process $Z\ \ $is\ a delayed regenerative process with
regeneration points $\{\overset{\ast }{T}_{n},n=0,1,...\}$.$\ $

\bigskip

We define the bivariate process $B=(Z,R)$, where for $t\geq 0\ \ R_{t}\ \ $%
is the release rate ($0$ or $M$)\ at time $t$. The process $B$ has as its
state space the pair of line segments

\bigskip 
\begin{equation*}
S=[(l,\lambda )\times \{0\}]\cup \lbrack (\tau ,V]\times \{M\}]\text{,}
\end{equation*}%
where $l\ \ $is the lower bound of the state space of the input process $I$.

\bigskip \smallskip

The penalty cost rate function is given by

\bigskip

\begin{equation*}
f(z)=\left\{ 
\begin{array}{cc}
g(z)\text{,} & (z,r)\in (l,\lambda )\times \{0\} \\ 
g^{\ast }(z)\text{,} & (z,r)\in (\tau ,V]\times \{M\}%
\end{array}%
\right. \ \ (3.2)
\end{equation*}%
where $g:(l,\lambda )\rightarrow R_{+}$ and $g^{\ast }:(\tau ,V]\rightarrow
R_{+}$ are bounded measurable function.

\bigskip

For $\alpha \in R_{+}$, let the $C_{\alpha }(x,\lambda ,0)\ $and $C_{\alpha
}(x,\tau ,M)$ be the expected discounted penalty costs during the interval $%
[0,\overset{\wedge }{T}_{0})$, and during the interval $[\overset{\wedge }{T}%
_{0},\overset{}{\overset{\ast }{T}}_{0})$ starting at $x$, respectively. It
follows that, for $x\in \lbrack \tau ,V]$

\begin{equation*}
C_{\alpha }(x,\lambda ,0)=E_{x}\int_{0}^{\overset{\wedge }{T}_{0}}e^{-\alpha
t}g(I_{t})dt\text{,\ \ \ }
\end{equation*}%
and for $x\in \lbrack \lambda ,V]$%
\begin{equation*}
C_{\alpha }(x,\tau ,M)=E_{x}\int_{0}^{\overset{}{\overset{\ast }{T}}%
_{0}}e^{-\alpha t}g^{\ast }(I_{t}^{\ast })dt\text{. \ \ }(3.3)
\end{equation*}

\bigskip

The functionals (3.3), $E_{x}[e^{-\alpha \overset{\wedge }{T}_{0}}]$, $E_{x}[%
\overset{\wedge }{T}_{0}]$, $E_{x}[e^{-\alpha \overset{}{\overset{\ast }{T}}%
_{0}}]$, $E_{x}[\overset{}{\overset{\ast }{T}}_{0}]$, which we aim to
evalaute, are needed to obtain the total discounted and the long-run average
costs associated with the $P_{\lambda ,\tau }^{M}\ $policy, discussed in
Section 4.

\bigskip \smallskip

For any $a\in R$, we define $T_{a}^{+}=\inf \{t\geq 0:\overset{}{I}_{t}\
\geq a\}$, $T_{a}^{-}=\inf \{t\geq 0:\overset{}{I}_{t}\ \leq a\}$, $\top
_{a}^{+}=\inf \{t\geq 0:\overset{\wedge }{I}_{t}\ \geq a\}$ and $\top
_{a}^{-}=\inf \{t\geq 0:\overset{\wedge }{I}_{t}\ \leq a\}$.\ We note that $%
\overset{\wedge }{T}_{0}=T_{\lambda }^{+}$\ almost everywhere.\ \ \ 

\smallskip

To derive $C_{\alpha }(x,\lambda ,0)$, $E_{x}[e^{-\alpha \overset{\wedge }{T}%
_{0}}]$, $E_{x}[\overset{\wedge }{T}_{0}]\ $we define the process obtained
by \ killing the process $I$\ \ at $\overset{\wedge }{T}_{0}$, as follows:

\begin{equation*}
X_{t}=\{I_{t},t<\overset{\wedge }{T}_{0}\}\text{.}\ (3.4)
\end{equation*}%
It is known that this killed process is a strong Markov process, with state
space $(l,\lambda )$.

\smallskip

For any Borel set $A\ \subset (l,\lambda )$, and $t\in R_{+}$,$\ $the
probability transition function of this process is given as follows

\begin{equation*}
P_{t}(x,A)=P_{x}(I_{t}\in A,t<\overset{\wedge }{T}_{0})
\end{equation*}%
and for each $\alpha \in R_{+}$ its $\alpha -$potential is defined as follows

\begin{equation*}
U^{\alpha }(x,A)=\int_{0}^{\infty }P_{t}(x,A)e^{-\alpha t}dt=E_{x}\int_{0}^{%
\overset{\wedge }{T}_{0}}e^{-\alpha t}\boldsymbol{I}_{\{I_{t}\in A\}}dt\text{%
. \ }(3.5)
\end{equation*}%
We note that for $x<\lambda $

\begin{equation*}
C_{\alpha }(x,\lambda ,0)=U^{\alpha }g(x)\text{. \ }\ (3.6)
\end{equation*}

\ The following Theorem will be used extensively throughout this paper.

\smallskip

\textbf{Theorem 1}. Let $S=\{S_{t},t\geq 0\}\ $ be a strong Markov process.
Define, $\mathcal{G}=\left\{ \sigma (S_{u}\ ,u\leq t)\right\} _{t\geq 0}$,%
{\LARGE \ }$\tau $\ to be any stopping time with respect to $\mathcal{G}\ $.
Let Y be the process obtained by killing the process S\ at time $\tau $,
denote the state space of this process by $E$, and let $U^{\alpha }\ $be\
its $\alpha -$potential. Then, for $x\in E$ 
\begin{equation*}
E_{x}[e^{-\alpha \tau }]=1-\alpha U_{{}}^{\alpha }\mathbf{I}_{E}(x)\text{.\
\ \ }(3.7)
\end{equation*}

\textbf{Proof}. From the definition of $U^{\alpha }$,$\ $\ for any bounded \
measurable function $f$\ \ whose domain is $E$, we have

\begin{eqnarray*}
U^{\alpha }f(x) &=&E_{x}[\int_{0}^{\tau }e^{-\alpha
t}f(S_{t})dt]=\int_{E}f(y)U^{\alpha }(x,dy)\cdot \\
&&
\end{eqnarray*}%
Taking $f\ \ $to be identically equal to one, we have

\begin{equation*}
\frac{1-E_{x}[e^{-\alpha \tau }]}{\alpha }=U_{{}}^{\alpha }\mathbf{I}_{E}(x)%
\text{.}
\end{equation*}%
The required result is immediate from the last equation above. \ \ $%
\blacksquare $ \ 

\bigskip

\bigskip\ \ \ \ \ \ \ \ \ \ \ \ \ \ \ \ \ \ \ \ \ \ \ \ \ \ 

First we consider the case where the input process$\ $is a spectrally
positive\ L\'{e}vy process.

\textbf{Proposition 1}. For $\alpha \geq 0$,\ $a\leq \lambda $ the \ $\alpha 
$-potential $\ $($\overset{(1)}{U^{\alpha }}$)$\ $of the process $I\ \ $%
killed at\ $T=\overset{\wedge }{T}_{0}\wedge $\ $T_{a}^{-}\ $is absolutely
continuous with respect to the Lebesgue measure on $[a,\lambda ]\ $and a
version of its density is given by\ 

\begin{equation*}
\overset{(1)}{u^{\alpha }}(x,y)=W^{(\alpha )}(\lambda -x)\frac{W^{(\alpha
)}(y-a)}{W^{(\alpha )}(\lambda -a)}-W^{(\alpha )}(y-x),\ \ x,y\in \lbrack
a,\lambda ]\text{.}\ \ \ \ \ \ (3.8)
\end{equation*}

\textbf{Proof}. For $A\subset \lbrack a,\lambda ]$

\begin{eqnarray*}
\overset{(1)}{U^{\alpha }}(x,A) &=&E_{x}\int_{0}^{T}e^{-\alpha t}\mathbf{I}%
_{\{I_{t}\in A\}}dt \\
&=&E_{-x}\int_{0}^{\top _{-\lambda }^{-}\wedge \top
_{-a}^{+}}e_{{}}^{-\alpha t}\mathbf{I}_{\{\overset{\wedge }{I}_{t}\in -A\}}dt
\\
&=&E_{\lambda -x}\int_{0}^{\top _{0}^{-}\wedge \top _{\lambda
-a}^{+}}e_{{}}^{-\alpha t}\mathbf{I}_{\{\overset{\wedge }{I}_{t}\in \lambda
-A\}}dt \\
&=&\int_{(\lambda -A)}{\Large [}W^{(\alpha )}(\lambda -x)\frac{W^{(\alpha
)}(\lambda -a-y)}{W^{(\alpha )}(\lambda -a)}-W^{(\alpha )}(y-x)]dy\text{,}
\end{eqnarray*}%
where the last equation follows from Theorem 8.7 of [7], this establishes
our assertion. $\ \ \ \ \ \ \ \ \ \ \ \blacksquare $\ \ 

\ \ \ \ \ \ \ \ \ \ \ \ \ \ \ \ \ \ \ \ 

\textbf{Corollary 1}. For $\alpha \geq 0$ the $\alpha $-potential $%
(U^{\alpha })\ $of the process $X$\ \ is absolutely continuous with respect
to the Lebesgue measure on $(-\infty ,\lambda ]\ $and a version of its
density is given by

\bigskip 
\begin{equation*}
u^{\alpha }(x,y)=W^{(\alpha )}(\lambda -x)e^{-\eta (\alpha )(\lambda
-y)}-W^{(\alpha )}(y-x),\ \ x,y\in (-\infty ,\lambda ]\text{. \ }(3.9)
\end{equation*}

\textbf{Proof}. The proof follows from (3.8) by letting $a\rightarrow
-\infty $ and since,\ for $\alpha \geq 0$, $W^{(\alpha )}(x)\sim \frac{%
e^{\eta (\alpha )x}}{\phi ^{^{\prime }}(\eta (\alpha )\ )}\ $as $%
x\rightarrow \infty $. \ \ \ \ $\ \blacksquare $\ $\ $

\bigskip

We are now in a position to find $E_{x}[e^{-\alpha \overset{\wedge }{T}%
_{0}}] $ and $E_{x}[\overset{\wedge }{T}_{0}]$.

\bigskip \textbf{Proposition 2}. (i) For $\alpha \geq 0$ and $x\leq \lambda
\ $we$\ $have

\begin{equation*}
E_{x}[e^{-\alpha \overset{\wedge }{T}_{0}}]=Z^{(\alpha )}(\lambda -x)-\frac{%
\alpha }{\eta (\alpha )}W^{(\alpha )}(\lambda -x)\text{. }(3.10)
\end{equation*}

\ \ \ \ \ \ \ \ \ \ \ \ (ii) For $x\leq \lambda \ $we$\ $have%
\begin{eqnarray*}
E_{x}[\overset{\wedge }{T}_{0}] &=&\frac{W(\lambda -x)}{\eta (0)}-\overset{-}%
{W}(\lambda -x)\text{, }\eta (0)>0\text{\ \ } \\
&=&\infty \ \ \ \ \ \ \ \ \ \ \ \ \ \ \ \ \ \ \ \ \ \ \ \ \ \ \ \text{, }%
\eta (0)=0\text{, \ }(3.11)
\end{eqnarray*}%
where for every $x\geq 0$,%
\begin{equation*}
\overset{-}{W}(x)=\int_{0}^{x}W(y)dy.\ \ 
\end{equation*}

\textbf{Proof}. We only prove (i), the proof of (ii)\ is easily obtained
from (i)\ and hence is omitted. Let $\overset{}{U^{\alpha }}\ $be as defined
in Corollary 1, then

\begin{eqnarray*}
E_{x}[e^{-\alpha \overset{\wedge }{T}_{0}}] &=&1-\alpha U_{{}}^{\alpha }%
\mathbf{I}_{(-\infty ,\lambda )}(x) \\
&=&1-\alpha \int_{-\infty }^{\lambda }{\large \{}W^{\alpha }(\lambda
-x)e^{-(\lambda -y)\eta (\alpha )}-W^{\alpha }(y-x){\Large \}}dy \\
&=&1+\alpha \int_{x}^{\lambda }W^{\alpha }(y-x)dy-\alpha W^{\alpha }(\lambda
-x)\int_{-\infty }^{\lambda }e^{-(\lambda -y)\eta (\alpha )}dy \\
&=&Z^{(\alpha )}(\lambda -x)-\frac{\alpha }{\eta (\alpha )}W^{(\alpha
)}(\lambda -x)\text{,}
\end{eqnarray*}%
where the first equation follows from (3.7), the second equation follows
from (3.9), the third equation follows since $W^{(\alpha )}(x)=0$, $x\,<0$,\
and the last equation follows from the definition of $Z^{(\alpha )}$. \ \ \
\ $\blacksquare \ \ $

\smallskip

For any Borel set $B\subset R_{+}\times R$, we let $M(B)$ be the Poisson
random measure counting the number of jumps of the process $I$ \ in $B$ with
L\'{e}vy measure $\nu $,$\ $where if $B=[0,t)\times A$, $A\subset R$, then $%
E[M(B)]=t\upsilon (A)$. \ \ \ 

\smallskip\ \ \ \ \ \ \ \ \ \ \ \ \ \ \ \ \ \ \ \ \ \ \ 

\textbf{Proposition 3}. Let $\overset{(1)}{u^{\alpha }}\ $be as given in
(3.8), and $x\leq \lambda \leq z$, then 
\begin{equation*}
E_{_{x}}[e^{-\alpha \overset{\wedge }{T}_{0}},I_{\overset{\wedge }{T}%
_{0}}\in dz,\ \overset{\wedge }{T}_{0}<T_{a}^{-}]=\int_{a}^{\lambda
}\upsilon (dz-y)\overset{(1)}{u^{\alpha }}(x,y)dy\ \ \ \ \ (3.12)
\end{equation*}

\textbf{Proof}. Let $T$ \ be as defined in Proposition 1. For $x<\lambda
,\alpha \geq 0$, $C\subset \lbrack \lambda ,\infty )\ $and $D\subset
(a,\lambda )$ we have 
\begin{eqnarray*}
E_{x}[e^{-\alpha \overset{\wedge }{T}_{0}},\text{ }I_{\overset{\wedge }{T}%
_{0}} &\in &C,I_{\overset{\wedge }{T}_{0}-}\in D,\overset{\wedge }{T}%
_{0}<T_{a}^{-}] \\
&=&E_{x}{\LARGE [}\int_{[0,\infty )\times (0,\infty )}e_{{}}^{-\alpha t}%
\mathbf{I}_{\{\overset{-}{I}_{t}-<\lambda ,\underline{I_{t-}}>a,I_{t}-\in
D\}}\mathbf{I}_{\{y\in C-I_{t}-\}}M(dt,dy){\Large ]} \\
&=&E_{x}{\LARGE [}\int_{[0,\infty )}e^{-\alpha t}\mathbf{I}_{\{\overset{-}{I}%
_{t-}<\lambda ,\underline{I_{t-}}>a\}}\mathbf{I}_{\{I_{t}\in D\}}\nu
(C-I_{t})dt{\Large ]} \\
&=&E_{x}{\LARGE [}\int_{[0,\infty )}e^{-\alpha t}\mathbf{I}_{\{t<T\}}\nu
(C-I_{t})\mathbf{I}_{\{I_{t}\in D\}}dt){\Large ]} \\
&=&E_{x}{\LARGE [}\int_{[0,\infty )\times D}e^{-\alpha t}\mathbf{I}%
_{\{t<T\}}\nu (C-y)\mathbf{I}_{\{I_{t}\in dy\}}dt{\Large ]} \\
&=&\int_{D}\nu (C-y)\overset{(1)}{u^{\alpha }}(x,y)dy,
\end{eqnarray*}%
where the second equation follows from the \textit{compensation formula\ }%
(Theorem 4.4. of [7])\textit{. }Our assertion is proved by taking $%
D=[a,\lambda ]$. \ \ \ $\underset{}{_{{}}\underset{}{}}_{{}}$\ \ \ \ \ $%
\blacksquare $ \ \ \ \ \ \ \ \ 

\ \ \ \ \ \ \ \ \ \ \ \ \ \ \ \ \ \ \ \ \ \ \ \ \ \ \ 

The following corollary gives a formula needed to compute the total
discounted cost.

\smallskip

\textbf{Corollary 2}. Let $u^{\alpha }\ $be\ as defined in (3.9). For $%
\alpha \geq 0$ and for $x\leq \lambda \leq z$,

\begin{equation*}
E_{_{x}}[e^{-\alpha \overset{\wedge }{T}_{0}},I_{\overset{\wedge }{T}%
_{0}}\in dz\ ]=\int_{-\infty }^{\lambda }\upsilon (dz-y)u^{\alpha }(x,y)dy%
\text{. }\ \ \ (3.13)
\end{equation*}

\textbf{Proof}. The proof follows immediately from (3.9) and (3.12) by
letting $a\rightarrow -\infty $.\ \ \ \ \ \ \ \ \ \ \ \ \ \ \ \ \ \ \ \ \ \
\ $\blacksquare $ \ 

\smallskip \smallskip \medskip

We now turn our attention to the case where the input process is a
spectrally positive L\'{e}vy process reflected at its infimum. In this case,
the killed process has state space $[0,\lambda )$. Let $\overset{(2)}{%
U^{\alpha }}$\ be the $\ \alpha $-$\ $potential of this process.

\ 

\textbf{Proposition 4}. For any $x,y\in \lbrack 0,\lambda )$,

\ \ \ \ \ \ \ \ \ \ \ \ \ \ \ \ \ \ \ \ \ \ \ \ \ \ \ \ \ \ \ \ \ \ \ \ \ \
\ \ \ \ \ \ \ \ \ \ \ \ \ \ \ \ \ \ \ \ \ \ \ \ \ \ \ \ \ \ \ \ \ \ \ \ \ \
\ \ \ \ \ \newline
\begin{equation*}
\overset{(2)}{U^{\alpha }}(x,dy)=\frac{W^{(\alpha )}(\lambda -x)W^{(\alpha
)}(dy)}{W_{+}^{(\alpha )^{\prime }}(\lambda )}-W^{(\alpha )}(y-x)dy\text{, \
\ \ \ }(3.14)
\end{equation*}%
where for $x,y\in \lbrack 0,\lambda )$, $W^{(\alpha )}(dy)=W^{(\alpha
)}(0)\delta _{0}(dy)+W_{+}^{(\alpha )^{\prime }}(y)dy$, and $\delta _{0}\ $%
is the delta measure in zero.

\ \ \ \ \textbf{Proof}. Note that for each $t\geq 0,$%
\begin{eqnarray*}
I_{t} &=&Y_{t}-\underline{Y_{t}}\ \ \ \ \ \ (3.15)\ \  \\
&=&\overset{-}{\overset{\wedge }{Y}}_{t}-\overset{\wedge }{Y}_{t}\text{,}
\end{eqnarray*}%
where the process $Y=\{Y_{t},t\geq 0\}\ $is a spectrally positive L\'{e}vy
process. The result follows from part (ii) \ of Theorem 1\ of [10], since
the process $\overset{}{\overset{\wedge }{Y}}\ $is a spectrally negative L%
\'{e}vy process.~\ \ \ $\blacksquare $

\bigskip

The following provides results parallel to (3.10) and (3.11), resectively.

\smallskip

\textbf{Proposition 5}. Assume that the input process is a spectrally
positive L\'{e}vy process reflected at its infimum. Then

(i) For $\alpha \geq 0$ and $x\leq \lambda \ $we$\ $have%
\begin{equation*}
E_{x}[e^{-\alpha \overset{\wedge }{T}_{0_{{}}}}]=Z^{(\alpha )}(\lambda
-x)-W^{(\alpha )}(\lambda -x)\frac{\alpha W^{(\alpha )}(\lambda )}{%
W_{+}^{(\alpha )^{\prime }}(\lambda )}.\ \text{ }(3.16)
\end{equation*}

(ii) For $x\leq \lambda \ $we$\ $have%
\begin{equation*}
E_{x}[\overset{\wedge }{T}_{0}]=W(\lambda -x)\frac{W(\lambda )}{%
W_{+}^{^{\prime }}(\lambda )}-\overset{-}{W}(\lambda -x)\text{. \ }(3.17)
\end{equation*}

\textbf{Proof}. The proof of part (i) follows from (3.7) and (3.14), in a
manner similar to the proof of (3.10). The proof of part (ii) follows from
part (i) by direct differentiation. $\ \ \ \ \ \ \ \blacksquare $

\bigskip

To find a formula analogous to (3.13), when the input is a spectrally
positive\ L\'{e}vy process reflected at its infimum, we first need few
definitions. Define

\bigskip 
\begin{eqnarray*}
l_{\alpha }(dz) &=&W^{(\alpha )}(\lambda -x)\int_{0}^{\lambda }W^{(\alpha
}(dy)\upsilon (dz-y)\ \ \ \ \ \ \ \ \ \ \ \ \ \ \ \ \ \  \\
&&-W_{+}^{(\alpha )^{\prime }}(\lambda )\int_{0}^{\lambda }dyW^{(\alpha
)}(y-x)\upsilon (dz-y)],z>\lambda .\ (3.18) \\
L_{\alpha }(z) &=&\int_{(z,\infty )}^{{}}l_{\alpha }(du).\text{\ \ \ \ \ \ \
\ \ \ \ \ \ \ \ \ \ \ \ \ \ \ \ \ \ \ \ \ \ \ \ \ \ \ \ \ \ \ \ \ \ \ \ \ \
\ }(3.19) \\
V_{\alpha }(\lambda ) &=&W_{+}^{(\alpha )^{\prime }}(\lambda )Z^{(\alpha
)}(\lambda -x)-\alpha W^{(\alpha )}(\lambda -x)W^{(\alpha )}(\lambda ).\text{%
\ \ \ }(3.20)\ 
\end{eqnarray*}

The following proposition gives the required formula.

\smallskip

\textbf{Proposition 6}. (i) For $\alpha \geq 0$ and for $x\leq \lambda <$ $z$%
,

\begin{equation*}
E_{_{x}}[e^{-\alpha \overset{\wedge }{T}_{0}},I_{\overset{\wedge }{T}%
_{0}}\in dz]=\frac{l_{\alpha }(dz)}{W_{+}^{(\alpha )^{\prime }}(\lambda )}\ 
\text{.\ }\ \ \ \ \ (3.21)
\end{equation*}

\ \ \ \ \ \ \ \ \ \ \ \ \ \ \ \ \ (ii) For $\alpha \geq 0$

\begin{equation*}
E_{_{x}}[e^{-\alpha \overset{\wedge }{T}_{0}},I_{\overset{\wedge }{T}%
_{0}}=\lambda ]=\frac{V_{\alpha }(\lambda )-L_{\alpha }(\lambda )}{%
W_{+}^{(\alpha )^{\prime }}(\lambda )}.\ \ (3.22)\ 
\end{equation*}

\smallskip

\textbf{Proof}. (i) Consider the spectrally positive L\'{e}vy process $%
Y=\{Y_{t},t\geq 0\}$, given in the proof of Proposition 4.$\ $For any $a\in
R $, we define $\digamma _{a}$ \ as the sigma algebra generated by $(Y_{s}\
,s\leq t)$, $\tau _{a}^{+}=\inf \{t\geq 0:Y_{t}\ \geq a\}$, $\tau
_{a}^{-}=\inf \{t\geq 0:Y_{t}\ \leq a\}$, $\sigma _{a}^{+}=\inf \{t\geq 0:%
\overset{\wedge }{Y}_{t}\ \geq a\}$,~and $\sigma _{a}^{-}=\inf \{t\geq 0:%
\overset{\wedge }{Y}_{t}\ \leq a\}$. From (3.15), for $x\geq 0$, $I_{0}\ =x\
\ $if and only if $Y_{0}=x$ \ if and only if $\overset{\wedge }{Y}_{0}=-x$ .
Furthermore, $\overset{\wedge }{T}_{0}=\tau _{\lambda }^{+}$ and $I_{\overset%
{\wedge }{T}_{0}}=Y_{\tau _{\lambda }^{+}}$ \ almost surely on $\{\tau
_{\lambda }^{+}<\tau _{0}^{-}\}$.\ Therefore

\begin{eqnarray*}
E_{_{x}}[e^{-\alpha \overset{\wedge }{T}_{0}},I_{\overset{\wedge }{T}_{0}}
&\in &dz]=E_{_{x}}[e^{-\alpha \overset{\wedge }{T}_{0}},I_{\overset{\wedge }{%
T}_{0}}\in dz,\tau _{\lambda }^{+}<\tau _{0}^{-}]+E_{_{x}}[e^{-\alpha 
\overset{\wedge }{T}_{0}},I_{\overset{\wedge }{T}_{0}}\in dz,\tau _{\lambda
}^{+}\geq \tau _{0}^{-}] \\
&=&E_{_{x}}[e^{-\alpha \tau _{\lambda }^{+}}\ ,Y_{\tau _{\lambda }^{+}}\in
dz,\tau _{\lambda }^{+}<\tau _{0}^{-}] \\
&&+E_{x}[e^{-\alpha \tau _{0}^{-}},\tau _{\lambda }^{+}\geq \tau
_{0}^{-}]\times E_{_{0}}[e^{-\alpha \overset{\wedge }{T}_{0}},I_{\overset{%
\wedge }{T}_{0}}\in dz] \\
&=&E_{_{x}}[e^{-\alpha \tau _{\lambda }^{+}}\ ,Y_{\tau _{\lambda }^{+}}\in
dz,\tau _{\lambda }^{+}<\tau _{0}^{-}] \\
&&+E_{-x}{\Large [}e^{-\alpha \sigma _{0}^{+}},\sigma _{-\lambda }^{-}\geq
\sigma _{0}^{+}{\Large ]}\times E_{_{0}}{\Large [}e^{-\alpha \overset{\wedge 
}{T}_{0}},I_{\overset{\wedge }{T}_{0}}\in dz{\large ]} \\
&=&E_{_{x}}[e^{-\alpha \tau _{\lambda }^{+}}\ ,Y_{\tau _{\lambda }^{+}}\in
dz,\tau _{\lambda }^{+}<\tau _{0}^{-}]\ \ \ \ \ \ \ \ \ \ \  \\
&&+E_{\lambda _{-x}}[e^{-\alpha \sigma _{\lambda }^{+}},\sigma
_{0}^{-}>\sigma _{\lambda }^{+}]\times E_{_{0}}[e^{-\alpha \overset{\wedge }{%
T}_{0}},I_{\overset{\wedge }{T}_{0}}\in dz]\text{{\large , \ \ \ \ \ \ \ \ \
\ \ \ \ \ \ \ \ \ \ \ \ \ \ \ \ \ \ \ \ \ \ \ \ \ \ }}
\end{eqnarray*}%
where the second equation follows from the first equation by conditioning on 
$\digamma _{\tau _{0}^{-}}\ $and then using the strong Markov property. The
third and fourth equations follow from the definitions of $\overset{}{%
\overset{\wedge }{Y}\text{, }}$ $\tau _{a}^{+}$,$\tau _{a}^{-}\ $,$\sigma
_{a}^{+}$,$\sigma _{a}^{-}$. \ 

Letting $a\rightarrow 0\ $in (3.8) and (3.12), we find that the first term
in the last equation above is equal to $\dint\limits_{0}^{\lambda }\nu (dz-y)%
{\Large [}W^{(\alpha )}(\lambda -x)\frac{W^{(\alpha )}(y)}{W^{(\alpha
)}(\lambda )}-W^{(\alpha )}(y-x){\large ]}dy$. The second term is equal to $%
\frac{W^{(\alpha )}(\lambda -x)}{W^{(\alpha )}(\lambda )}\ $(see (8.8) of
[7]) and the third term is equal to $\frac{h_{\alpha }(dz)}{W_{+}^{(\alpha
)^{\prime }}(\lambda )}$ (this follows from Theorem 4.1 of [11] by letting
the $\beta ,\gamma \rightarrow 0$).

Our assertion is satisfied by replacing each of the three terms in the last
equation\ by the corresponding value indicated above and after some
algebraic manipulations, which we omit.

\ \ \ \ \ \ \ \ \ (ii) The proof is immediate from (3.16) and (3.21). \ \ \
\ \ \ \ \ \ \ \ \ \ \ \ \ \ $\ \blacksquare $ \ 

\ 

\bigskip Now we turn our attention to computing $C_{\alpha }(x,\tau ,M)$, $%
E_{x}[\exp (-\alpha \overset{\ast }{T}_{0})]$, and $E_{x}[\overset{\ast }{T}%
_{0}]$, when $x\in \lbrack \lambda ,V]$.

Let ${\large \eta }_{\tau }=\inf \{t\geq 0:\overset{\ast }{I}_{t}\leq \tau
\} $ and,$\ $for each $t\geq 0$,
\begin{equation*}
\overset{\ast }{X}_{t}=\{\overset{\ast }{I}_{t},t<{\large \eta }_{\tau }\}%
\text{.}\ \ \ (3.23)
\end{equation*}%
Note that, the state space of the process $\overset{\ast }{X}\ $is the
interval $(\tau ,V]$, and let $\overset{\ast }{U}^{\alpha }$be its $\alpha $%
-potential. Starting at any $x\in \lbrack \lambda ,V]$, ${\large \eta }%
_{\tau }=\overset{\ast }{T}_{0}\ $almost everywhere, furthermore the sample
paths of a spectrally positive L\'{e}vy process and a spectrally positive L%
\'{e}vy process reflected at its infimum behave the same way until they
reach level $\tau $, thus $\overset{\ast }{X}\ $behaves the same way in both
cases. It follows that, for each $x\in \lbrack \lambda ,V]$, 
\begin{equation*}
C_{\alpha }(x,\tau ,M)=\overset{\ast }{U}^{\alpha }g^{\ast }(x)\text{.}\ \ \
(3.24)
\end{equation*}

Denote the process $I-M$\ \ by $N$, note that this process is a spectrally
positive L\'{e}vy process with the L\'{e}vy exponent $\phi _{M}(\theta
)=\phi (\theta )+\theta M$,$\ \theta \bigskip \geq 0$. $\ $We denote its $%
\alpha -$scale and adjoint $\alpha -$scale functions by $W_{M}^{(\alpha )}$
and $Z_{M}^{(\alpha )}$, respectively.

\medskip

\textbf{Theorem 2}. For $\alpha \geq 0,$ $\overset{\ast }{U}^{\alpha }$is
absolutely continuous with respect to the Lebsegue measure on $(\tau ,V]$,
and a version of its density is given by

\begin{equation*}
\overset{\ast }{u}^{\alpha }(x,y)=\frac{Z_{M}^{(\alpha )}(V-x)W_{M}^{(\alpha
)}(y-\tau )}{Z_{M}^{(\alpha )}(V-\tau )}-W_{M}^{(\alpha )}(y-x)\ \text{,}\ \
\ x,y\in (\tau ,V]\text{. \ }(3.25)
\end{equation*}

\textbf{Proof.} For each $t\geq 0$, we define $B_{t}=N_{t}-V$.\ $\ \ $For\
any $b\in R$, we define $\sigma _{b}^{-}=\inf \{t\geq 0:B_{t}-\overset{-}{%
B_{t}}_{{}}<b\}\ $and$\ \gamma _{b}^{+}=\inf \{t\geq 0:\overset{\boldsymbol{%
\wedge }}{B_{t}}-\underline{\overset{\boldsymbol{\wedge }}{B_{t}}}>b\}.$ For
any Borel set $A\subseteq (\tau ,V]$ and $x\in (\tau ,V]$ we have 
\begin{eqnarray*}
P_{x}\{\overset{\ast }{X}_{t} &\in &A\}=P_{x}\{\overset{\ast }{I}_{t}\in A,t<%
{\large \eta }_{\tau }\} \\
&=&P_{x}\{N_{t}-\underset{s\leq t}{\sup }((N_{s}-V)\vee 0)\in A,t<{\large %
\eta }_{\tau }\} \\
&=&P_{x-V}\{B_{t}-\overset{\boldsymbol{-}}{B\ }_{t}\in A-V,t<\sigma _{\tau
-V}^{-}\} \\
&=&P_{V-x}\{\overset{\boldsymbol{\wedge }}{B}_{t}-\underline{\overset{%
\boldsymbol{\symbol{94}}}{B}_{t}}\ \in V-A,t<\gamma _{V-\tau }^{+}\}
\end{eqnarray*}%
Using Theorem 1 (i) of [10], the result follows. \bigskip\ \ \ $\blacksquare 
$

The following theorem gives Laplace transform of the distribution of the
stopping time $\overset{\ast }{T}_{0}\ $and $E_{x}[\overset{\ast }{T}_{0}]$
when $x\in \lbrack \lambda ,V]$.

\bigskip

\textbf{Theorem 3.} (i) Let $x\in \lbrack \lambda ,V]\ $and $\alpha \in
R_{+} $, then

\begin{equation*}
E_{x}[e^{-\alpha \overset{\ast }{T}_{0}}]=\frac{Z_{M}^{(\alpha )}(V-x)}{%
Z_{M}^{(\alpha )}(V-\tau )}\text{. \ }(3.26)
\end{equation*}

\ \ \ \ \ \ \ \ \ \ \ \ \ \ \ \ \ \ \ \ (ii)\ \ \ For $x\in \lbrack \lambda
,V]$%
\begin{equation*}
E_{x}[\overset{\ast }{T}_{0}]=\overset{-}{W_{M}}(V-\tau )-\overset{-}{W_{M}}%
(V-x)\text{,\ \ \ }(3.27)
\end{equation*}%
where, $\overset{-}{W_{M}}(x)=\ \dint\limits_{0}^{x}\overset{}{W_{M}}(y)dy$.

\textbf{Proof}. \ We only prove (i), the proof of (ii) follows easily from
(i) and is omitted. For $x\in \lbrack \lambda ,V]$, we have

\begin{eqnarray*}
E_{x}[e^{-\alpha \overset{\ast }{T}_{0}}] &=&1-\alpha \overset{\ast }{U}%
^{\alpha }\mathbf{I}_{(\tau ,V]}(x) \\
&=&1-\alpha \int_{\tau }^{V}\overset{\ast }{u}^{\alpha }(x,dy) \\
&=&1-\alpha \int_{\tau }^{V}{\Large [}\frac{Z_{M}^{(\alpha )}(V-x)W^{\alpha
}(y-\tau )}{Z_{M}^{(\alpha )}(V-\tau )}-W_{M}^{(\alpha )}(y-x){\Large ]}dy\ 
\\
&=&1-\alpha {\Huge [}\frac{Z_{M}^{(\alpha )}(V-x)}{Z_{M}^{(\alpha )}(V-\tau )%
}{\Large \{}\frac{Z^{\alpha }(V-\tau )-1}{\alpha }{\Large \}}-{\large \{}%
\frac{Z_{M}^{(\alpha )}(V-x)-1}{\alpha }{\large \}}{\Huge ]}
\end{eqnarray*}%
\begin{eqnarray*}
&=&\frac{Z_{M}^{(\alpha )}(V-x)}{Z_{M}^{(\alpha )}(V-\tau )}-Z_{M}^{(\alpha
)}(V-x)+Z_{M}^{(\alpha )}(V-x) \\
&=&\frac{Z_{M}^{(\alpha )}(V-x)}{Z_{M}^{(\alpha )}(V-\tau )}\text{,}
\end{eqnarray*}%
where the third equation follows from (3.25), the fourth equation follows
from the definition of the function $Z_{M}^{(\alpha )}\ $and the fifth
equation follows the fourth equation after obvious manipulations. \ $%
\blacksquare $

\bigskip

\textbf{Remark 1}. When $V=\infty $, for $\alpha \geq 0\ $we let $\eta
_{M}(\alpha )=\sup \{\theta :\phi (\theta )-\theta M=\alpha \}$. Since $%
Z_{M}^{(\alpha )}(y)=O(e^{\eta _{M}(\alpha )y})$ as $y\rightarrow \infty ,\ $%
then we have

\begin{eqnarray*}
E_{x}[\overset{\ast }{T}_{0}] &=&\frac{(x-\tau )}{\eta _{M}(\alpha
)^{^{\prime }}(0)} \\
&=&\frac{(x-\tau )}{M-E(I_{1})}\ ,\ \text{if }M>E(I_{1}) \\
&=&\infty \ \ \ \ \ \ \ \ \ \ \ \ \ \ \ \ ,\ \text{if\ }M\leq E(I_{1})\text{.%
}
\end{eqnarray*}%
This is consistent with the well known fact about the busy period of the
M/G/1 queuing system.

\bigskip The following gives $E_{x}[\exp (-\alpha \overset{\ast }{T}_{0})]$,
when $x<\lambda $, a result that is needed to compute the total discounted
cost.

\bigskip \textbf{Theorem 4}. Assume that the input process is a spectrally
positive L\'{e}vy\textit{\ } process. For $z>\lambda $, we define

\bigskip 
\begin{equation*}
h_{\alpha }(x,dz)=\int_{-\infty }^{\lambda }u^{\alpha }(x,y)\upsilon (dz-y)dy%
\text{,}
\end{equation*}%
where $u^{\alpha }(x,y)$ is defined in (3.9).

Then, for $\alpha \geq 0,x<\lambda $

\begin{equation*}
E_{x}[e^{-\alpha \overset{\ast }{T}_{0}}]=\frac{1}{Z_{M}^{(\alpha )}(V-\tau )%
}[\int_{\lambda }^{V}Z_{M}^{(\alpha )}(V-z)h_{\alpha
}(x,dz)+\int_{V}^{\infty }h_{\alpha }(x,dz)]\ \ \ \ (3.28)\ 
\end{equation*}

\bigskip \textbf{Proof}. \ We write\ \ \ \ \ \ \ \ \ \ \ \ \ \ \ \ \ \ \ \ \
\ \ \ \ \ \ \ \ \ \ \ \ \ \ \ \ \ \ \ \ \ \ \ \ \ \ \ \ \ \ \ \ \ \ \ \ \ \
\ \ \ \ \ \ \ \ \ \ \ \ \ \ \ \ \ \ \ \ \ \ \ \ \ \ \ \ \ \ \ \ \ \ \ \ \ \
\ \ \ \ \ \ \ \ \ \ \ \ \ \ \ \ \ \ \ \ \ \ \ \ \ \ \ \ \ \ \ \ \ \ \ \ \ \
\ \ \ \ \ \ \ \ \ \ \ \ \ \ \ \ \ \ \ \ \ \ \ \ \ \ \ \ \ \ \ \ \ \ \ \ \ \
\ \ \ \ \ \ \ \ \ \ \ \ \ \ \ \ \ \ \ \ \ 

\begin{eqnarray*}
E_{x}[e^{-\alpha \overset{\ast }{T}_{0}}] &=&E_{x}[e^{-\alpha \overset{%
\wedge }{T}_{0}-\alpha (\overset{\ast }{T}_{0}-\overset{\wedge }{T}_{0})}] \\
&=&E_{x}{\Large [}E_{x}[e^{-\alpha \overset{\wedge }{T}_{0}-\alpha (\overset{%
\ast }{T}_{0}-\overset{\wedge }{T}_{0})}\mid \sigma (\overset{\wedge }{T}%
_{0},I_{\overset{\wedge }{T}_{0}}){\Large ]} \\
&=&E_{x}{\Large [}e^{-\alpha \overset{\wedge }{T}_{0}\ \ \ }E_{(I_{\overset{%
\wedge }{T}_{0}}\wedge V)}[e^{-\alpha \overset{\ast }{T}_{0}}]{\Large ]} \\
&=&\frac{1}{Z_{M}^{(\alpha )}(V-\tau )}E_{x}[e^{-\alpha \overset{\wedge }{T}%
_{0}\ \ }Z_{M}^{(\alpha )}(V-(I_{\overset{\wedge }{T}_{0}}\wedge V))] \\
&=&\frac{1}{Z_{M}^{(\alpha )}(V-\tau )}[\int_{\lambda }^{V}Z_{M}^{(\alpha
)}(V-z)h_{\alpha }(x,dz)+\int_{V}^{\infty }h_{\alpha }(x,dz)]\text{,}
\end{eqnarray*}%
where the third equation follows since, given $\overset{\wedge }{T}_{0}$ and 
$I_{\overset{\wedge }{T}_{0}}$, $\overset{\ast }{T}_{0}-\overset{\wedge }{T}%
_{0}\ $is equal to $\overset{\ast }{T}_{0}\ $almost everywhere. The fourth
equation follows from (3.26). The last equation follows from (3.13), the
fact that $\ Z_{M}^{(\alpha )}(0)=1$, and the definition of $h_{\alpha
}(x,dz)$.$\ $ \ \ $\blacksquare $

\bigskip

The following theorem gives a result analogous to (3.28) when the input
process is a spectrally positive L\'{e}vy process reflected at its infimum.

\bigskip \textbf{Theorem 5}. Assume that the input process is a spectrally
positive L\'{e}vy process\ reflected at its infimum. For $z\geq \lambda $,
let $l_{\alpha }(dz)$, $L_{\alpha }(z)$, and $V_{\alpha }(\lambda )$\ be as
defined in (3.18), (3.19) and (3.20), respectively. Define

\begin{equation*}
g_{\alpha }(x,dz)=\left\{ 
\begin{array}{c}
=\frac{l_{\alpha }(dz)}{W_{+}^{(\alpha )^{\prime }}(\lambda )}\ \text{,\ }%
z>\lambda \\ 
=\frac{V_{\alpha }(\lambda )-L_{\alpha }(\lambda )}{W_{+}^{(\alpha )^{\prime
}}(\lambda )}\delta _{\lambda }(dz)\text{.}%
\end{array}%
\right.
\end{equation*}

\ \ Then, for $\alpha \geq 0,x<\lambda $

\begin{equation*}
E_{x}[e^{-\alpha \overset{\ast }{T}_{0}}]=\frac{1}{Z_{M}^{(\alpha )}(\lambda
-\tau )}[\int_{\lambda }^{V}Z_{M}^{(\alpha )}(\lambda -z)g_{\alpha
}(x,dz)+\int_{V}^{\infty }g_{\alpha }(x,dz)]\text{. \ \ \ \ }(3.29)\ \ \ 
\end{equation*}

\textbf{Proof}. The proof follows in a manner similar to the proof of
(3.28), using (3.21), (3.22) and (3.26). \ \ \ \ \ $\blacksquare $\ \ \ \ \
\ \ \ 

\bigskip \medskip

\medskip \bigskip

{\LARGE 4.} {\Large The expected total discounted and long-run average costs}

\medskip

Consider a finite dam controlled by a $P_{\lambda ,\tau }^{M}$ policy as
described in Section 1. Assume that the input process, $I$, is spectrally
positive L\'{e}vy, and define $\alpha $ to be the discount factor. For $x\in
\lbrack \tau ,V]$, we let $C_{x}^{\alpha }(\lambda ,\tau )$, and $C(\lambda
,\tau )$ be the expected total discounted cost and long-run average cost,
respectively, given $I_{0}=0$. Furthermore, we define $C^{\alpha }(x)$ as
the expected discounted cost during the interval $[0,\overset{\ast }{T}_{0})$%
, given the initial water content is equal to $x$.

\bigskip

Modifying (3.1) of [1], it follows that for $x\in \lbrack \tau ,V]$,%
\begin{equation*}
C_{x}^{\alpha }(\lambda ,\tau )=C^{\alpha }(x)+\frac{E_{x}[\exp (-\alpha 
\overset{\ast }{T}_{0})]C^{\alpha }(\tau )}{1-E_{\tau }[\exp (-\alpha 
\overset{\ast }{T}_{0})]}\text{.}\ (4.1)
\end{equation*}

\bigskip

From the definition of the $P_{\lambda ,\tau }^{M}\;$policy, it follows that
for $\lambda <x<V$

\begin{equation*}
C^{\alpha }(x)=M\{K_{1}-RE_{_{x}}\int_{0}^{\overset{\ast }{T}_{0}}e^{-\alpha
t}dt\}+\;C_{\alpha }^{{}}(x,\tau ,M)\text{, }\ \ (4.2)
\end{equation*}%
and for $x\in \lbrack \tau ,\lambda ]$

\begin{eqnarray*}
C_{\alpha }(x) &=&M\{K_{2}+\ K_{1}\ E_{x}[e^{-\alpha \overset{\wedge }{T}%
_{0}}]-\frac{R}{\alpha }\{E_{_{x}}[e^{-\alpha \overset{\wedge }{T}_{0}}]\
-E_{_{x}}[e^{-\alpha \overset{\ast }{T}_{0}}]\ \} \\
&&+C_{\alpha }(x,\lambda ,0)+E_{_{x}}[e^{-\alpha \overset{\wedge }{T}%
_{0}}C_{\alpha }((I_{\overset{\wedge }{T}_{0}}\wedge V),\tau ,M)]\text{,}\ \
\ \ \ \ (4.3)\ \ \ 
\end{eqnarray*}%
where $C_{\alpha }(x,\lambda ,0)$, and $C_{\alpha }(x,\tau ,M)\ $ are given
in (3.6) and (3.24), respectively. Using (3.9), (3.10), (3.13), (3.25),
(3.26) and (3.28) we obtain $C^{\alpha }(x)$. Finally, the expected total
discounted cost can be determined explicitly by substituting (4.2), (4.3),
(3.26) and (3.28) into (4.1).

\smallskip

\ To determine the long-run average cost using a given $P_{\lambda ,\tau
}^{M}$ policy, we proceed as follows. Let $C(\lambda ,\tau )\ $denote the
long-run average cost, and $\ $define $C^{0}(x)$ as the expected
non-discounted cost during the interval $[0,\overset{\ast }{T}_{0})$, given
the initial water content is equal to $x$, $x\in \lbrack \tau ,V]$. It
follows that

\bigskip\ 
\begin{equation*}
C(\lambda ,\tau )=\frac{C^{0}(\tau )}{E_{\tau }[\overset{\ast }{T}_{0}]}.\
(4.4)
\end{equation*}%
From the strong Markov property we have

\begin{equation*}
E_{\tau }[\overset{\ast }{T}_{0}]=E_{\tau }[\overset{\wedge }{T}%
_{0}]+E_{\tau }{\Large [}E_{_{(I_{\overset{\wedge }{T}_{0}}\wedge V)}}%
{\large [}\overset{\ast }{T}_{0}{\large ]}{\Large ]}\text{. \ }\ (4.5)
\end{equation*}%
Furthermore

\begin{equation*}
C^{0}(\tau )=M{\LARGE \{}K-R{\Large (}E_{\tau }[\overset{\ast }{T}_{0}\
]-E_{\tau }[\overset{\wedge }{T}_{0}]{\Large )}{\LARGE \}}+C_{0}(\tau
,\lambda ,0)+E_{_{\tau }}[C_{0}((I_{\overset{\wedge }{T}_{\lambda }}\wedge
V),\tau ,M)]\text{,\ }\ \ (4.6)
\end{equation*}%
where $K=$\ $K_{1}+\ K_{2}\ $. Letting $\alpha =0$ in (3.13) and
substituting the result, along with\ (3.11) and (3.27) into (4.5) we obtain $%
E_{\tau }[\overset{\ast }{T}_{0}]$. Using (3.6), (3.9), (3.11), (3.13),
(3.24), (3.25) and (4.5) we obtain (4.6). Substituting (4.5) and (4.6) into
(4.4) the long-run average cost is determined.

\bigskip

The corresponding results for the spectrally positive L\'{e}vy reflected at
its infimum input follow similarly.\ \ \ \ \ \ \ 

\bigskip\ \ 

\medskip \newpage

\ \ \ \ \ \ \ \ \ \ \ \ \ \ \ \ \ \ \ \ \ \ \ \ \ \ \ \ \ \ \ \ \ \ \ \ 

\bigskip \noindent {\LARGE 5. Special Cases}

\medskip {\LARGE \ }

In this section we consider the cases where the input process is a
spectrally positive L\'{e}vy \ of bounded variation, Brownian motion
reflected at its infimum and Wiener process. For the first case, we extend
the results of [6], we also simplify some of their results. For the second
case, we obtain results similar to those of [4] and [8]. In the third case
we obtain the results of [12].

\smallskip \smallskip

\textbf{Case 1}{\LARGE . }Assume that the input is a spectrally positive L%
\'{e}vy process of bounded variation with L\'{e}vy exponent described in
(2.3), reflected at its infimum. Let \noindent\ $\mu =\int_{0}^{\infty
}x\upsilon (dx)$ and assume that\ $\mu <\infty $. For every $x\in R_{+}$, we
define the probability density function $f(x)=\frac{\upsilon ([x,\infty ))}{%
\mu }$. We have $\int_{0}^{\infty }(1-e^{-\theta x})\upsilon (dx)=$ $\theta
\mu \int_{0}^{\infty }\ e^{-\theta x}f(x)dx$.\ Define $\rho =\frac{\mu }{%
\varsigma }$, $F(x)\ $as\ the distribution function corresponding to $f\ $.
Assume that $\rho <1$, it follows that, $\frac{1}{\phi (\theta )}=\frac{1}{%
\varsigma }\int_{0}^{\infty }e^{-\theta x}dx\sum\limits_{n=0}^{\infty }\rho
^{n}F^{(n)}(x)$. Therefore, the $0-$scale function is given is given as
follows

\begin{equation*}
W(x)=\frac{1}{\varsigma }\sum\limits_{n=0}^{\infty }\rho ^{n}F^{(n)}(x)\text{%
. \ }(5.1)\text{\ }
\end{equation*}%
For $\alpha >0,W^{(\alpha )}$ is computed using (2.5)\ and (5.1).

\bigskip

Define $\overset{\ast }{\varsigma }=\varsigma +M$, and $\overset{\ast }{\rho 
}=\frac{\mu }{\overset{\ast }{\varsigma }}$. Let $W_{M}^{(\alpha )}\ $be as
defined in the paragrah proceeding Theorem 2, and denote $W_{M}^{(0)}$ by $%
W_{M}^{{}}\ $,\ using an argument similar to the one above we have

\begin{equation*}
W_{M}^{{}}(x)=\frac{1}{\overset{\ast }{\varsigma }}\sum\limits_{n=0}^{\infty
}\overset{}{\overset{}{\overset{\ast }{\rho }^{n}}^{{}}}F^{(n)}(x)\text{. \ }%
(5.2)\text{\ \ }
\end{equation*}%
$\ $

\bigskip

It follows that, for a spectrally positive L\'{e}vy process of bounded
variation, $W(0)=0$. Thus, the $\alpha $-potential $\overset{(2)}{U^{\alpha }%
}\ $(given in (3.14) is absolutely continuous. From (4.4) (4.5), and (4.6),
the long-run average cost is determined once for $x\in \lbrack \lambda ,V]$,$%
\ E_{\tau }[\overset{\wedge }{T}_{0}]$, $E_{x}[\overset{\ast }{T}_{0}]$, $%
\overset{(2)}{U^{0}}$, $\overset{\ast }{u}^{0}$, and the didtribution of $I_{%
\overset{\wedge }{T}_{0}}$are computed. Using (5.1)\ and (3.17) we compute $%
E_{\tau }[\overset{\wedge }{T}_{0}]$, and using (5.2) and (3.27) $E_{x}[%
\overset{\ast }{T}_{0}]$ is determined for $x\in \lbrack \lambda ,V]$.
Furthermore, $\overset{(2)}{U^{0}}$ is computed using (3.14) and (5.1).\
From (3.25) it follows that, for $x,y\in (\tau <V]$,$\
u^{0}(x,y)=W_{M}^{{}}(y-\tau )-W_{M}^{{}}(y-x)$, which is determined using
(5.2). The distribution of $I_{\overset{\wedge }{T}_{0}}$is given by letting 
$\alpha \rightarrow 0$ in (3.21) and (3.22), and using (3.18), (3.19),
(3.20), and (5.1).

\bigskip

The corresponding results for the total discounted cost follow similarly.

\bigskip

\textbf{Remark 2}.

\medskip

(a) Bae\textit{\ et al} \ [6] obtain the long-run average cost, when the
input process is a compound Poisson process with a negative drift. In this
case, $\upsilon (dx)=\lambda G(dx)$,$\,\ $where $\lambda >0$ and $G\ $\ is a
distribution function of a positive random variable $[0,\infty )$,
describing the size of each jump of the compound Poisson process. $\ $In
this case, $f(x)=\frac{\overset{-}{G}(x)\ }{m}$ and $\rho =\frac{\lambda m}{%
\varsigma },\ $where $\overset{-}{G}=1-G$\ \ and $m$ $=\int_{0}^{\infty }%
\overset{-}{G}(x)dx$, which is assumed to be finite. We note that their
entities $w(x)$ and $E[L^{a}(\lambda ,\tau )]$ given in p.521 and p.523 of
[6], respectively, are nothing but our $C_{0}(x,\lambda ,0)$ and $E_{\tau }[%
\overset{\wedge }{T}_{0}]$, respectively. Furthermore, for $x\in \lbrack
\tau ,V]$, their functions $u(x)\ $and $E[T_{\tau }^{M}\ _{{}}(x)]\ $given
in page 524 are our $C_{0}(x,\tau ,M)\ $and $E_{x}[\overset{\ast }{T}_{0}]$,
respectively. The distribution of $L(\tau )\ $(the overshoot) given on page
525 follows in an obvious manner from the distribution of \ $I_{\overset{%
\wedge }{T}_{0}}$. The formulas for computing $C_{0}(x,\lambda ,0)$, $%
E_{\tau }[\overset{\wedge }{T}_{0}]$, $C_{0}(x,\tau ,M)$, $E_{x}[\overset{%
\ast }{T}_{0}]$, and the distribution of $I_{\overset{\wedge }{T}_{0}}\ $,
follow from the corresponding results obtained in Case 1. We note that our
formulas for computing $w(x)$ and $E[L^{a}(\lambda ,\tau )]$ are identical
to those of [6], while our formulas for $u(x)$,\ $E[T_{\tau }^{M}\ _{{}}(x)]$%
, and the distribution of $L(\tau )\ $are simpler than theirs$.$

(b) Assume that the input process is a gamma process with negative drift.
The L\'{e}vy measure is given by $\upsilon (dx)=a\frac{e^{-bx}}{x}dx${\LARGE %
, }$a,b>0.\ $In this case, $E(I_{1})=\varsigma -\frac{a}{b}$, which is
assumed to be nonnegative and $\rho =\frac{a}{\varsigma b}$\ $<1$.$\ $It
follows that $f(x)=b\int_{x}^{\infty }${\LARGE \ }$\frac{e^{-by}{\LARGE \ }}{%
y}dy$, the \ right hand side is denoted by $E_{1}(x)\ $in p. 227 of [3], \
Direct integrations yield, $F(x)=(1-e^{-bx})+xf(x)$.$\ \ \ $\ \ \ \ 

(c) Assume that the input process is an inverse Gaussian process with a
negative drift, and with L\'{e}vy measure is given by $\upsilon (dx)=\frac{1%
}{\sigma \sqrt{2\pi x^{3}}}e^{-xc^{2}/2\sigma ^{2}}${\LARGE , }$\sigma ,c>0$%
. It follows that $E(I_{1})=\varsigma -\frac{1}{c}$, which is assumed to be
greater than zero. In this case $\rho =$ $\frac{1}{c\varsigma }<1$, $%
f(x)=c\int_{x}^{\infty }{\LARGE \ }\upsilon (dy)$, and $F(x)=\func{erf}(c%
\sqrt{y/2\sigma ^{2}})+xf(x)$.\ \ \ \ \ \ \ \ \ \ \ \ \ \ \ \ \ \ \ \ \ \ \
\ \ \ \ \ \ 

(d) Since, for all $n\geq 0$, $F^{(n)}(x)\leq \lbrack F(x)]^{n}$, then for
all $x$, $W(x)\leq \frac{1}{\varsigma -\mu F(x)}$,\ if $\rho <1$.$\ \ \ \ $

$\ \ \ \ \ \ \ \ \ \ \ \ \ \ \ \ \ \ \ \ \ \ \ \ \ \ \ \ \ \ \ \ \ \ \ \ \ \
\ \ \ \ \ \ \ \ \ \ \ \ \ \ \ \ \ \ \ \ \ \ \ \ \ \ \ \ \ \ \ \ \ $

\textbf{Case 2}{\LARGE . }Assume that the input process is a Brownian motion
with drift term $\mu \in R$, variance term $\sigma ^{2}$, reflected at its
infimum. From (3.6), (3.24), and (4.1)-(4.6), the total discounted and
long-run average costs are determined once $\ E_{x}[\overset{}{e^{-\alpha 
\overset{\wedge }{T}_{0}}}]$, $\ E_{\lambda }[\overset{}{e^{-\overset{\ast }{%
\alpha T}_{0}}}]$, $E_{\tau }[\overset{\wedge }{T}_{0}]$, $E_{\lambda }[%
\overset{\ast }{T}_{0}]$, $\overset{(2)}{U^{\alpha }}$, and $\overset{\ast }{%
u}^{\alpha }\ $are computed.

\bigskip

\ In this case, $I_{\overset{\wedge }{T}_{0}}=\lambda \ \ $almost
everywhere, the L\'{e}vy measure $\nu =0$, and from (2.1) we have, for $%
\theta \geq 0$, $\phi (\theta )=-\mu \theta +\frac{\theta ^{2}\sigma ^{2}}{2}
$. It follows that, for $\alpha \geq 0$, $\eta (\alpha )=\frac{\sqrt{2\alpha
\sigma ^{2}+\mu ^{2}}+\mu }{\sigma ^{2}}$. Let $\delta =\sqrt{2\alpha \sigma
^{2}+\mu ^{2}}$, we have, $W^{(\alpha )}(x)=\frac{2}{\delta }e^{\mu x/\sigma
^{2}}\sinh (\frac{x\delta }{\sigma ^{2}})\ $and $\ Z^{(\alpha )}(x)=e^{\mu
x/\sigma ^{2}}\left( \cosh (\frac{x\delta }{\sigma ^{2}})-\frac{\mu }{\delta 
}\sinh (\frac{x\delta }{\sigma ^{2}})\right) $.$\ \ $We note that $W^{\alpha
}(x)\ $is differentiable, and $W^{(\alpha )^{\prime }}(x)=\frac{\mu }{\sigma
^{2}}W^{(\alpha )}(x)+\frac{2}{\sigma ^{2}}e^{\mu x/\sigma ^{2}}\cosh (\frac{%
x\delta }{\sigma ^{2}})$; hence $\frac{W^{(\alpha )}(\lambda )}{W^{(\alpha
)^{\prime }}(\lambda )}=\left( \frac{\sigma ^{2}}{\mu +\delta \coth (\frac{%
\lambda \delta }{\sigma ^{2}})}\right) $. Substituting the values of $%
Z^{(\alpha )}(\lambda -x),W^{(\alpha )}(\lambda -x)$ and $\frac{W^{(\alpha
)}(\lambda )}{W^{(\alpha )^{\prime }}(\lambda )}\ $in (3.16), we have, for $%
\alpha \geq 0,x\leq \lambda $

\begin{equation*}
E_{x}[\overset{}{e^{-\alpha \overset{\wedge }{T}_{0}}}]=e^{\mu (\lambda -x)}%
\left[ \cosh \left( \frac{(\lambda -x)\delta }{\sigma ^{2}}\right) -\frac{1}{%
\delta }\sinh \left( \frac{(\lambda -x)\delta }{\sigma ^{2}}\right) \left(
\mu +\frac{2\alpha \sigma ^{2}}{\mu +\delta \coth (\frac{\lambda \delta }{%
\sigma ^{2}})}\right) \right] \text{, \ }(5.3)
\end{equation*}%
a simpler and more explicit formula than (4.6) of [4].

\bigskip We note that, $W^{(\alpha )}(0)=0$, $\overset{(2)}{U^{\alpha }}$\ $%
\ $is absolutely continuous \ with respect to the Lebesgue measure on $%
[0,\lambda )$. Substituting the values of $W^{(\alpha )}(\lambda -x)$,\ $%
W^{\alpha ^{\prime }}(\lambda )$, $W^{(\alpha )}(y-x)$, and $W^{\alpha
^{\prime }}(y)\ $in (3.14) we get a version of the the density of $\overset{%
(2)}{U^{\alpha }}$.

\smallskip

Let$\ \overset{\ast }{\mu }^{{}}=\mu -M$,$\ \overset{\ast }{\delta }=\sqrt{%
2\alpha \sigma ^{2}+\overset{\ast }{\mu }^{2}}$, \ it follows that $%
W_{M}^{(\alpha )}(x)=\frac{2}{\delta }e^{\overset{\ast }{\mu }x/\sigma
^{2}}\sinh (\frac{x\overset{\ast }{\delta }}{\sigma ^{2}})$, and $%
Z_{M}^{(\alpha )}(x)=e^{\overset{\ast }{\mu }x/\sigma ^{2}}\left( \cosh (%
\frac{x\overset{\ast }{\delta }}{\sigma ^{2}})-\frac{\mu }{\delta }\sinh (%
\frac{x\overset{\ast }{\delta }}{\sigma ^{2}})\right) $.\ Hence, $\overset{%
\ast }{u}^{\alpha }\ $is computed using (3.25). \ \ Since $I_{\overset{%
\wedge }{T}_{0}}=\lambda \ \ $almost everywhere,

\begin{equation*}
E_{x}[e^{-\alpha \overset{\ast }{T}_{0}}]=E_{x}[e^{-\alpha \overset{\wedge }{%
T}_{0}}]E_{\lambda }[\overset{}{e^{-\overset{}{\alpha }\overset{\ast }{T}%
_{0}}}]\text{, \ }(5.4)
\end{equation*}%
where $E_{x}[e^{-\alpha \overset{\wedge }{T}_{0}}]\ $is given in (5.3) and $%
E_{\lambda }[\overset{}{e^{-\overset{\ast }{\alpha T}_{0}}}]=\frac{%
Z_{M}^{(\alpha )}(V-\lambda )}{Z_{M}^{(\alpha )}(V-\tau )}$, which follows
from (3.26).$\ \ \ \ $

$\ \ \ $To compute $E_{\tau }[\overset{\wedge }{T}_{0}]$, $E_{\lambda }[%
\overset{\ast }{T}_{0}]$, we first assume that $\mu \#0$. In this case, for $%
x\geq 0,$ $W(x)=\frac{e^{2\mu x/\sigma ^{2}}-1}{\mu },W^{^{\prime }}(x)=%
\frac{2e^{2\mu x/\sigma ^{2}}}{\sigma ^{2}}$ and $\overset{-}{W}(x)=\frac{%
\sigma ^{2}}{2\mu ^{2}}(e^{2\mu x/\sigma ^{2}}-1)-\frac{x}{\mu }$.$\ $
Substituting the values of $W(\lambda -x)$, $\overset{-}{W}(\lambda
-x),W(\lambda )\ $and $W^{^{\prime }}(\lambda )\ $in (3.17) we have,

\begin{equation*}
E_{\tau }[\overset{\wedge }{T}_{0}]=\frac{\lambda -\tau }{\mu }+\frac{\sigma
^{2}}{2\mu ^{2}}\left[ e^{-2\mu \lambda /\sigma ^{2}}-e^{-2\mu \tau /\sigma
^{2}}\right] .\ (5.5)
\end{equation*}%
\ Since, $I_{\overset{\wedge }{T}_{0}}=\lambda \ \ $almost everywhere, from
(4.5)\ we have$\ \ \ \ \ \ \ \ \ \ \ \ \ \ \ \ \ \ \ \ \ \ \ \ \ \ \ \ \ \ \
\ \ \ \ \ \ \ \ \ \ \ \ \ \ \ \ \ \ \ \ \ \ \ \ \ \ \ \ \ \ \ \ \ \ \ \ \ \
\ \ \ \ \ \ \ \ \ \ \ \ \ \ \ \ \ \ \ \ \ \ \ \ \ \ \ \ \ \ \ \ \ \ \ \ \ \
\ \ \ \ \ \ \ \ \ \ \ \ \ \ \ \ \ \ \ \ \ \ \ \ \ \ \ \ \ \ \ \ \ \ \ \ \ \
\ \ \ \ \ \ \ \ \ \ \ \ \ \ \ \ \ \ \ \ \ \ \ \ \ \ \ \ \ \ \ \ \ \ \ \ \ \
\ \ \ \ \ \ \ \ \ \ \ \ \ \ \ \ \ \ \ \ \ \ \ \ \ \ \ \ \ \ \ \ \ \ \ \ \ \
\ \ \ \ \ \ \ \ \ \ \ \ \ \ \ \ \ \ \ \ \ \ \ \ \ \ \ \ \ \ \ \ \ \ \ \ \ \
\ \ \ \ \ \ \ \ \ \ \ \ \ \ \ \ \ \ \ \ \ \ \ \ \ \ \ \ \ \ \ \ \ \ \ \ \ \ $

\bigskip 
\begin{equation*}
E_{\tau }[\overset{\ast }{T}_{0}]=E_{\tau }[\overset{\wedge }{T}%
_{0}]+E_{\lambda }[\overset{\ast }{T}_{0}].\ \ (5.6)
\end{equation*}%
We note that $\overset{}{\overset{-}{W}_{M}}$ is computed from $\overset{-}{W%
}\ $above by replacing the term $\mu \ \ $by $\overset{\ast }{\mu }$ defined
in the preceding paragraph. Let$\ $ $\lambda ^{\ast }=V-\lambda $, and $\tau
^{\ast }=V-\tau $, substituting the values $\overset{-}{W}_{M}(V-\tau )\ $%
and $\overset{-}{W}_{M}(V-\lambda )\ $in (3.27) we have

\begin{equation*}
E_{\lambda }[\overset{\ast }{T}_{0}]=\frac{\lambda ^{\ast }-\tau ^{\ast }}{%
\mu ^{\ast }}+\frac{\sigma ^{2}}{2\overset{\ast }{\mu }^{2}}\left[ e^{2%
\overset{\ast }{\mu }\tau ^{\ast }/\sigma ^{2}}-e^{2\mu \lambda ^{\ast
}/\sigma ^{2}}\right] \text{.\ }(5.7)
\end{equation*}

\bigskip If $\mu =0$, $\ $then $\overset{}{W}(x)=\frac{2x}{\sigma ^{2}}$ and 
$\overset{-}{W}(x)=\frac{x^{2}}{\sigma ^{2}}$. From (3.17) it follows that

\begin{equation*}
E_{\tau }[\overset{\wedge }{T}_{0}]=\frac{\lambda ^{2}-x^{2}}{\sigma ^{2}}%
\text{. \ }(5.8)
\end{equation*}%
It is easily shown that, \ $\overset{}{W_{M}}(x)=\frac{1}{M}{\Large (}$ $%
1-e^{-2xM/\sigma ^{2\ }})$, using (3.27) we have

\begin{equation*}
E_{\lambda }[\overset{\ast }{T}_{0}]=\frac{\tau ^{\ast }-\lambda ^{\ast }}{M}%
+\frac{\sigma ^{2}}{2M^{2}}\left[ e^{-2M\tau ^{\ast }/\sigma
^{2}}-e^{-2M\lambda ^{\ast }/\sigma ^{2}}\right] .\ (5.9)
\end{equation*}

\bigskip \smallskip \medskip We note that our (5.5) is consistent with (4.9)
of [4], while (5.8) is identical to the corresponding equation given in p.
298 of the same reference.

\smallskip \medskip \smallskip

\noindent \textbf{Case 3}{\LARGE .} Assume that the input process is a
Brownian motion with drift term $\mu \ >0$\ and variance parameter $\sigma
^{2}$. It follows that $\eta (0)=\frac{2\mu }{\sigma ^{2}}$. Substituting
the values of $W^{(\alpha )}(x),Z^{(\alpha )}(x)$, given in Case 2, in
(3.10) we have, for $x\leq \lambda $, $E_{x}[e^{-\alpha \overset{\wedge }{T}%
_{0}}]=\exp \left( (\delta -\mu )(x-\lambda )\right) $. Substituting $\frac{1%
}{\mu }(e^{2\mu x/\sigma ^{2}}-1)$ and $\frac{\sigma ^{2}}{2\mu ^{2}}%
(e^{2\mu x/\sigma ^{2}}-1)-\frac{x}{\mu }$ for $W(x)\ $and $\overset{-}{W}(x)
$, respectively, in (3.11), we have, for $x\leq \lambda $, $E_{x}[\overset{%
\wedge }{T}_{0}]=\frac{\lambda -x}{\mu }$. These results are consistent with
the results of Zuckerman [12], p.423. The values of $E_{\lambda }[\overset{}{%
e^{-\overset{}{\alpha }\overset{\ast }{T}_{0}}}]$ and $E_{\lambda }[\overset{%
\ast }{T}_{0}]\ $are the same whether the input process is spectrally
positive L\'{e}vy or spectrally positive L\'{e}vy reflected at its infimum.
The computations of the total discounted and long-run average costs can be
obtained using (3.9), (3.13), (3.25), (4.1)-(4.6)  and the values of $E_{x}[%
\overset{}{e^{-\alpha \overset{\wedge }{T}_{0}}}]$, $E_{\lambda }[\overset{}{%
e^{-\overset{}{\alpha }\overset{\ast }{T}_{0}}}]$, $E_{x}[\overset{\wedge }{T%
}_{0}]$,$\ E_{\lambda }[\overset{\ast }{T}_{0}]$, in manners similar to
those discussed in Case 2. 

\bigskip \bigskip

{\LARGE References}

\smallskip \medskip

[1] Abdel-Hameed, M. (2000). Optimal control of a dam using $P_{\lambda
,\tau }^{M}\ $policies and penalty cost when the input process is a compound
Poisson process with positive drift. \textit{J.Appl.Prob. }\textbf{37},
508-416.

[2] Abramov, V.M. (2007). Optimal control of a large dam. \textit{%
J.Appl.Prob.} \textbf{44}, 249-258.

[3] Abramowitz, M. and Stegun, I. A. (1964). \textit{Handbook of
Mathematical Functions}. Dover, New York.

[4] Attia, F. (1987). The control of a finite dam with penalty cost
function: Wiener process input. \textit{Stochastic Processes and Their
Applications} \textbf{25}, 289-299.

[5] Bae, J, Kim, S. and Lee, E.Y. (2002). A\ $P_{\lambda }^{M}$\ policy for
an M/G/1 queueing system. \textit{Appl. Math. Modelling} \ \textbf{26, }%
929-939.

[6] Bae, J, Kim, S. and Lee, E.Y. (2003). Average cost under $P_{\lambda
,\tau }^{M}\ $- policy in a finite dam with compound Poisson input. \textit{%
J.Appl.Prob. }\textbf{40}, 519-526.

[7] Kyprianou, A. E. (2006). \textit{Introductory Lecture Notes on
Fluctuations of L\'{e}vy Processes with Applications}. Springer Verlag.

[8] Lam, Y. and Lou, J.H. (1987). Optimal control of a finite dam: Wiener
process input. \textit{J.Appl.Prob.} \textbf{35}, 482-488.

[9] Lee, E.Y. and Ahn, S.K. (1998). $P_{\lambda }^{M}$\ policy for a dam
with input formed by a compound Poisson process. \textit{J.Appl.Prob.} 
\textbf{35}, 482-488.

[10] Pistorius, M.R. (2004). On exit and ergodicity of the sepctrally
one-sided L\'{e}vy process refelected at its infimum\textit{.}\ \textit{J.
Theor. Proab}. \textbf{17}, 183-220

[11] Zhou, X. W. (2007). Exit problems for sepctrally negative L\'{e}vy
processes refelected at either the supremum or the infimum. \textit{%
J.Appl.Prob.} \textbf{44}, 1012-1030.

[12] Zuckerman, D. (1977). Two-stage output procedure for a finite dam. 
\textit{J.Appl.Prob.} \textbf{14}, 421-425.

\bigskip\ \ \ \ 

{\LARGE \ }

\end{document}